\documentclass[11pt]{amsart}
\usepackage{}
\usepackage{amssymb}
\theoremstyle{plain}
\newtheorem{Thm}{Theorem}

\newtheorem{Cor}[Thm]{Corollary}

\begin{document}

\title[classical semilinear parabolic equation]
{ Boundary value problem for a classical semilinear parabolic
equation}

\author{Li Ma}

\address{ Department of mathematics \\
Henan Normal university \\
Xinxiang, 453007 \\
China}

\email{nuslma@gmail.com}

\thanks{ The research is partially supported by the National Natural Science
Foundation of China 10631020 and SRFDP 20090002110019}

\begin{abstract}
In this paper, we study the boundary value problem of the classical
semilinear parabolic equations
$$
u_t-\Delta u=|u|^{p-1}u, \ \ in \ \ \Omega\times (0,T)
$$
and $u=0$ on the boundary $\partial\Omega\times [0,T)$ and $u=\phi$
at $t=0$, where $\Omega\subset R^n$ is a compact $C^1$ domain,
$1<p\leq p_S$ is a fixed constant, and $\phi\in C^2_0(\Omega)$ is a
given smooth function. Introducing new idea, we show that there are
two sets $\tilde{W}$ and $\tilde{Z}$ such that for $\phi\in W$,
there is a global positive solution $u(t)\in \tilde{W}$ with $h^1$
omega limit $\{0\}$ and for $\phi\in \tilde{Z}$, the solution blows
up at finite time.


\end{abstract}

 \maketitle

\section{Introduction}
In this paper, we study the Dirichlet boundary value problem of the classical
semilinear parabolic equation
\begin{equation}\label{main-eq}
u_t-\Delta u=|u|^{p-1}u, \ \ \ in \ \ \Omega\times (0,T)
\end{equation}
with $u=0$ on the boundary $\partial\Omega\times [0,T)$ and $u=\phi$
at $t=0$, where $T>0$, $\Omega\subset R^n$ is a compact $C^1$
domain, $p>1$ is a fixed constant, and $\phi\in C^2_0(\Omega)$ is a
given smooth function. Assume that $p\leq p_S=\frac{n+2}{n-2}$ for
$n\geq 3$ and $p<\infty$ for $n=1,2$. By the standard theory we know
that there is a local time positive solution to (\ref{main-eq}).
With the help of Nehari functional, one may find the threshold of
the initial datum such that the solution either exists globally or
blows up in finite time. More interesting results about
(\ref{main-eq}) can be found in the recent work \cite{CDW}. Since
the equation (\ref{main-eq}) is a model problem, it deserves to have
more understanding. Introducing new idea, we show in this paper that
there are two new sets $\tilde{W}$ and $\tilde{Z}$ such that for
$\phi\in \tilde{W}$, there is a global positive solution in
$\tilde{W}$ with the $H^1$ omega limit $0$ and for $\phi\in
\tilde{Z}$, the solution blows up at finite time. We may extend the
method used in this paper to treat Neumann boundary value problem of
semilinear parabolic equation with negative power in \cite{MW}. To
define the invariant set $\tilde{Z}$, we shall use the fact that the
cones
$$C_+=\{u\in C_0^1(\Omega); u\geq 0, \ \ u\not=0\}$$ and $$C_-=\{u\in C_0^1(\Omega); u\leq 0, \ \ u\not=0\}$$ are invariant sets of (\ref{main-eq}). This fact can be proved by applying the maximum principle.

We now recall the standard way to construct the invariant sets for (\ref{main-eq}).
Formally, (\ref{main-eq}) has a Lyapunov functional; namely,
$$
J(u)=\int \frac{1}{2}|\nabla u|^2-\frac{1}{1+p}u^{p+1}.
$$
Here and after, we use $\int$ to denote the integration over $\Omega$.
 In fact, we may consider (\ref{main-eq}) as the negative
L2-gradient flow of the functional $J(\cdot)$. That is, abstractly,
(\ref{main-eq}) can be written as $$ u_t=-J'(u).
$$
Hence, we have
$$
\frac{d}{dt}J(u(t))=<J'(u),u_t>=-|u_t|^2_2=-|J'(u)|_{L^2}^2.
$$

 Let $f(u)=u^{p}$ and its primitive
$$
F(u)=\frac{u^{p+1}}{p+1}.
$$
 Introduce the working space
$$
\Sigma=\{u\in H^1_0; u\not=0, \int F(u)<\infty\}.
$$
The condition $\int F(u)<\infty$ is always true by using the Sobolev inequality.

 Define on $\Sigma$, the functional
 $$
M(u)=\frac{1}{2}\int_{\Omega}|u|^2
$$
and the Nehari functional
$$
I(u)=\int |\nabla u|^2-uf(u)=\int |\nabla u|^2-|u|^{p+1}.
$$
Note that these two functionals are well-defined on $\Sigma$.

 Along the flow (\ref{main-eq}) we can see that
\begin{equation}\label{mass}
\frac{d}{dt}M(u)=\int uu_t=-I(u).
\end{equation}

Let
$$
d=\inf\{J(u); u\in\Sigma; I(u)=0\}.
$$
Define
$$
W=\{u\in \Sigma; J(u)<d,
I(u)>0\}\bigcup\{0\}
$$ and
$$
Z=\{u\in \Sigma; J(u)<d,
I(u)<0\}.
$$
The classical result says that $W$ and $Z$ are invariant sets of (\ref{main-eq}); furthermore, for $1<p<p_S$ and for any initial data $\phi\in W$, the solution exists globally; for $1<p\leq p_S$ and for any initial data $\phi\in Z$, the solution blows up at finite time. One may see \cite{F} for more results and references.

We now introduce new functionals.
For $\lambda \in \mathbf{R}_+$, define
$$
E_{\lambda}(u)=J(u)+\lambda M(u).
$$
Then along the flow (\ref{main-eq}), we have
\begin{equation}\label{decrease}
\frac{d}{dt}E_{\lambda}(u)=-|J'(u)|^2_2-\lambda I(u).
\end{equation}
From this, it is clear that for $\lambda\geq 0$, we have
$$
\frac{d}{dt}E_{\lambda}(u)>0
$$
except $|J'(u)|=0$.

Introduce
$$
d_{\lambda}=\inf\{E_{\lambda}(u); u\in\Sigma; I(u)=0\}.
$$
As in the case for the quantity $d$, we can give it the
mountain-pass characterization.

Assume it is finite at this moment. Define
$$
W_{\lambda}=\{u\in \Sigma; E_{\lambda}(u)<d_{\lambda},
I(u)>0\}\bigcup\{0\}.
$$
For convenient we set $W_0=W$. Arguing as in $W$, one can see that $W_{\lambda}$ with $\lambda>0$ is non-empty.

Then by (\ref{decrease}) and the standard argument we know that for
$\lambda\geq 0$, $W_{\lambda}$ is a invariant set of the flow
(\ref{main-eq}).

One of our main results for (\ref{main-eq}) is to show the the following
conclusion.

\begin{Thm}\label{wei} Fix any power $1<p<p_S$, we have for $\lambda>0$ that

 (1). $d_{\lambda}$ is finite, and $d_\lambda>d$ for $\lambda>0$;

 (2). for $\phi\in W_{\lambda}$ with
$\lambda\geq 0$, the flow exists globally and its omega limit is
$\{0\}$. Hence $$\tilde{W}:=\bigcup_{\lambda\geq 0}W_{\lambda}$$ is
invariant set of (\ref{main-eq}).
\end{Thm}
We remark that since $d_\lambda>d$, we know that the set $W_\lambda$
is different from the set $W$.

To find the set for blow-up solutions to (\ref{main-eq}), we need to
use the comparison argument. We shall restrict the initial data
being positive. Let $\delta\geq 0$. Consider the boundary value
problem of the following semilinear parabolic equation
\begin{equation}\label{main-delta}
v_t-\Delta v+\delta v=v^{p}, \ \  u>0, \ \ in \ \ \Omega\times (0,T)
\end{equation}
with $u=0$ on the boundary $\partial\Omega\times [0,T)$ and $u=\phi$
at $t=0$, where $T:=T_{max}(\phi)>0$ is the maximal existence time
of the solution $v(t)$. Define on $\Sigma_+=\Sigma\bigcap C_+$,
$$
J_\delta(v)=J(v)+\delta M(v),
$$
$$
I_\delta(v)=I(v)+2\delta M(v),
$$
and on the set where $\{I_\delta(v)=0\}$
$$
E^{\delta}(v)=J_\delta(v)=(\frac{1}{2}-\frac{1}{p+1})\int |u|^{p+1}.
$$
Define
$$
d_\delta=\inf \{E^\delta(v); v\in \Sigma_+, I_\delta(u)=0\}.
$$

For $\epsilon>0$,
$$
d_{\delta, \epsilon}=\inf\{J_\delta(u); u\in\Sigma_+;
I_\delta(u)=\epsilon\}
$$
and
$$
Z_{\delta}=\{u\in \Sigma_+; J_\delta(u)<d_{\delta}, I_\delta(u)<0\}.
$$
Clearly, $Z_{\delta}$ is non-empty and it is a
invariant set of the flow (\ref{main-delta}). We remark that one may make similar construction on $\Sigma_-=\Sigma\bigcap C_-$.

\begin{Thm}\label{wei1} Fix $1<p\leq p_S$.
(1). For $\phi\in Z_{\delta}$ and $\phi\geq v_\delta$, the flow $(v(t))$ to (\ref{main-delta}) blows
up in finite time.

(2). Let $u(t)$ be the flow to (\ref{main-eq}) with the initial data $\phi$ as (1) above. Then $u(t)\geq v(t)$ and $u(t)$ blows up
at some $t<\infty$.

\end{Thm}

As a consequence of Theorem \ref{wei1}, we have
\begin{Cor}
Set $\tilde{Z}=\bigcup_{\delta\geq 0}Z_\delta$. Then for any $\phi\in \tilde{Z}$, the solution for (\ref{main-eq}) blows up at finite time.
\end{Cor}

The results above will be proved
in next section.
\section{Global solution and finite time blow-up solution}\label{sect1}

\textbf{We now prove Theorem \ref{wei}}.

(1). The finiteness of $d_{\lambda}$ can be obtained in the similar
way as in \cite{F}.  Since $1<p<p_S$, we know that
$d_\lambda$ can also be achieved by some function $u_\lambda$(see \cite{C}
\cite{M}, or \cite{S}).
 By this we know that $d_\lambda$ is
different from $d$ for $\lambda>0$. Hence,
we have $d_{\lambda}>d$ for $\lambda>0$.

(2). Since $I(\phi)>0$, we have $I(u(t))>0$ for all $t\in [0,T)$.
For otherwise, for some $t>0$, $I(u(t))=0$. Using the definition of
$d_\lambda$, we have $E_{\lambda}(u(t))\geq d_{\lambda}$. This is a
contradiction to the fact that
$$
\frac{d}{dt}E_{\lambda}(u(t))<0, \ \ and  \ \
E_{\lambda}(u(t))<E_{\lambda}(\phi)<d_\lambda.
$$

Using (\ref{mass}), we know that $M(u(t))<M(\phi)$. With the help of
the condition $E_{\lambda}(u(t))<d$ and $1<p<p_S$, we know that $u(t)\in H^1$ is
uniformly bounded and bounding constant depends only on $d$, $p$,
$|\Omega|$, and $M(\phi)$.

The $H^1$ omega limit at $t=\infty$ can be determined below. It is a
classical fact (\cite{F}) that the $H^1$ omega limit set $\omega(\phi)$
consists of classical equilibria. If $v\in \omega(\phi)$, we have
$I(v)=0$. If $v$ is nontrivial, we have
$$
E_\lambda(v)\geq d_\lambda.
$$
Impossible. Hence $v=0$, that is, $\omega(\phi)=\{0\}$.

This completes the proof of Theorem \ref{wei}.

The remaining part of this section we give the proof of Theorem \ref{wei1}.
\begin{proof}(Proof of Theorem \ref{wei1}).
Introduce
$$
A=\inf\{\frac{|\nabla u|^2_2+\delta|u|_2^2}{|u|_{p+1}^2};\ \  u\in
H_0^1(\Omega),\ \  u\not=0\}.
$$
Then it is easy to see that
$d_\delta=\frac{p-1}{2(p+1)}A^{(p+1)/(p-1)}$ (see \cite{C} \cite{M},
or \cite{S}). Assume that $0\not=v\in H_0^1(\Omega)$ such that
$I_\delta(v)=-\epsilon$. Then
\begin{equation}\label{star2}
E^{\delta}(v)=\frac{p-1}{2(p+1)}\int (|\nabla v|^2+\delta v^2)-\frac{\epsilon}{p+1}.
\end{equation}
Using the definition of $A$ we have
$$
\int (|\nabla v|^2+\delta v^2)\leq \int |v|^{p+1}\leq A^{-\frac{p+1}{2}}(\int (|\nabla v|^2+\delta v^2))^{\frac{p+1}{2}}.
$$
Hence,
$$
\int (|\nabla v|^2+\delta v^2)\geq A^{(p+1)/(p-1)}.
$$
Combining this with (\ref{star2}) we have
\begin{equation}\label{star3}
d_{\delta,\epsilon}\geq d_\delta-\frac{\epsilon}{p+1}.
\end{equation}

We now prove (1) in the statement of theorem \ref{wei1}.

(1). Take $\epsilon>0$ such that
$$
\epsilon<\min (-I_\delta(\phi),d_\delta-J_{\delta}(\phi)).
$$
Then using (\ref{decrease}) and (\ref{star3}) we know that
$$
J_{\delta}(v(t))\leq J_{\delta}(\phi)<d_{\epsilon}
$$
for $t\in [0,T)$. Since $I_\delta(\phi)<-\epsilon$, by using the definition
of $d_{\delta,\epsilon}$ and the continuity, we know that
$$
I_\delta(v(t))<-\epsilon.
$$
Note that $$
I_\delta(v)=2J_{\delta}(v)-(1-\frac{2}{p+1})\int |v|^{p+1}.
$$
Assume that $T=T_{max}>0$ be the maximal time of the flow $(v(t))$. Assume that $T=\infty$.
On one hand, using similar formula to (\ref{mass}) we have
$$
\frac{1}{2}\frac{d}{dt}\int v^2=-I_\delta(v)\geq \epsilon>0,
$$and then $$
\int v^2\geq \int \phi^2+2\epsilon t\to \infty.
$$
That is, $M(v(t))\to \infty$ as $t\to\infty$.

On the other hand,
$$
\frac{1}{2}\frac{d}{dt}\int v^2=-I_\delta(v)\geq -2d_\epsilon+(1-\frac{2}{p+1})\int |v|^{p+1}.
$$
Then we have
$$
\frac{d}{dt}M(v(t))\geq -2d_\epsilon+C(p,|\Omega|)M(v(t))^{\frac{p+1}{2}}
$$
for some uniform constant $C(p,|\Omega|)>0$. Then using $M(v(t))\to\infty$, we know that there exists
$T_1>0$ such that for any $t>T_1$,
$$
\frac{d}{dt}M(v(t))\geq \frac{1}{2}C(p,|\Omega|)M(v(t))^{\frac{p+1}{2}}.
$$
However, this implies that $T<\infty$. A contradiction. Hence $T<\infty$ and $M(v(t))\to\infty$ as $t\to T$.

We shall prove (2) in the statement of theorem \ref{wei1} by using
the comparison lemma. (2). Let $T_{max}<\infty$ be the blow-up time
of the flow $(v(t))$. Recall that $v(t)>0$ for $t\in (0,T_{max})$.
Let $w(t)=u(t)-v(t)$, $t<T_{\max}$. Then $w(t)$ is bounded in any
finite time before the blowing up time of the solution $u(t)$. Note
that
\begin{equation}\label{star}
w_t-\Delta w=p\xi^{p+1}w+\delta v.
\end{equation}
Recall that $w(0)=0$ and $w(t)|_{\partial \Omega}=0$. Let $w_-(t)$ be the negative part of $w(t)$.
Multiplying both sides of (\ref{star}) by $w_-(t)$ and integrating over $\Omega$ by $w_-(t)$, we get
$$
\frac{d}{dt}\int |w_-(t)|^2=-\int |\nabla w_-(t)|^2+p\int \xi^{-p-1} |w_-(t)|^2+\delta\int vw_-(t).
$$
We remark that the last term is non-positive. Then we have
$$
\frac{d}{dt}\int |w_-(t)|^2\leq C\int |w_-(t)|^2.
$$
By the Gronwall inequality we know that $\int|w_-(t)|^2=0$ for any $t>0$. Hence we have
$u(t)\geq v(t)$ and then $$\int u(t)^2\geq \int v(t)^2\to \infty$$ as $t\to T_{max}<\infty$.

\end{proof}

\end{document}